\documentclass[11pt,a4paper]{article}

\usepackage{amsmath,amsfonts,amssymb,amsthm}
\usepackage{bbm,mathrsfs, bm, mathtools}

\title{Random Gale diagrams and neighborly polytopes\\ in high dimensions}

\author{Rolf Schneider}
\date{}
\newcommand{\R}{{\mathbb R}}

\newcommand{\bP}{{\mathbb P}}

\newcommand{\N}{{\mathbb N}}

\newcommand{\D}{{\rm d}}

\newcommand{\bE}{{\mathbb E}\,}

\newtheorem{theorem}{Theorem}%[section]
\newtheorem{lemma}{Lemma}%[section]
\begin{document}
\maketitle

\begin{abstract}
Taking up a suggestion of David Gale from 1956, we generate sets of combinatorially isomorphic polytopes by choosing their Gale diagrams at random. We find that in high dimensions, and under suitable assumptions on the growth of the involved parameters, the obtained polytopes have strong neighborliness properties, with high probability.
\\[1mm]
\end{abstract}

\section{Introduction}\label{sec1}

The purpose of the following is, roughly speaking, to introduce a new class of random polytopes, which have strong neighborliness properties in high dimensions, with overwhelming probability. The main idea is to generate not directly the polytopes by a random procedure, but their Gale diagrams.

Let us first recall that a convex polytope $P$ in Euclidean space $\R^d$ is called $k$-neighborly if any $k$ or fewer vertices of $P$ are neighbors, which means that their convex hull is a face of $P$. A $\lfloor d/2\rfloor$-neighborly polytope is called neighborly. It is known that a $k$-neighborly $d$-polytope with $k>\lfloor d/2\rfloor$ must be a simplex, and that in dimensions $d\ge 4$ there are neighborly polytopes with any number of vertices. It appears that neighborliness properties (though not under this name) of certain polytopes were first noted by Carath\'{e}odory \cite{Car11}, and that the proper investigation of neighborly polytopes began with the work of Gale \cite{Gal56}. For more information on neighborly polytopes (for example, their important role in the Upper bound theorem), we refer to the books by Gr\"unbaum \cite{Gru03}, McMullen and Shephard \cite{McMS71}, Ziegler \cite{Zie95}, Matou\v{s}ek \cite{Mat02}.  

There is a widespread impression that there are ``many'' neighborly polytopes. For example, this is supported by the proof, due to Shemer \cite{She82}, that the number of combinatorial types of neighborly $(2m)$-polytopes with $2m+p$ vertices grows superexponentially as $p\to\infty$ ($m\ge 2$ fixed) and as $m\to\infty$ ($p\ge 4$ fixed).

Gale concluded his seminal paper \cite{Gal56} with the following `speculative remark', as he called it. ``It would seem that the likelihood of getting polyhedra every $m$ of whose vertices are neighbors increases rapidly with the dimension of the space.'' After briefly discussing the construction of a special $11$-dimensional polytope with 24 vertices via choosing in a random fashion the points of, as we would say nowadays, a Gale transform, he continued: ``Thus we might guess that finding every pair of points of $P$ neighbors would be the expected rather than the exceptional case. The probability notions hinted at here can be made precise and investigations in this direction would be of interest in further study of the problem.''

Surprisingly, it seems that this stimulus has not provoked a reaction so far. Of course, as remarked in \cite[p. 129b]{Gru03}, `the question ``What is the probability that a random (simplicial) polytope is neighborly?'' is not well-posed, since the answer heavily depends on the model of random polytopes.' The arguably most prominent model in geometric probability, the convex hull of (many) independent uniform random points in a given convex body, is not likely to exhibit strong neighborliness properties. On the other hand, a carefully chosen distribution and a suitably coordinated growth of dimension and number of random points may well show the desired effect. The following result, obtained by Donoho and Tanner \cite{DT05}, may serve as an example.

Let $T^{n-1}\subset\R^n$ denote the standard regular simplex in $\R^n$. Let $\Pi T^{n-1}$ be the image of $T^{n-1}$ under a  random orthogonal projection $\Pi$ from $\R^n$ to $\R^d$ (where $d<n$) with uniform distribution (that is, defined by a normalized Haar measure). Since the random polytope $\Pi T^{n-1}$, which is almost surely simplicial, has at most $n$ vertices, it can have at most $\binom{n}{k+1}$ faces of dimension $k$. Donoho and Tanner \cite{DT05} found an increasing function $\rho_N:(0,1)\to(0,1)$ (see \cite[Fig. 1]{DT05} for a picture of its graph) such that the following theorem holds. Here $f_k(P)$ is the number of $k$-dimensional faces of a polytope $P$, and $\bP$ denotes probability.

\vspace{2mm}

\noindent{\bf Theorem.} {\em If $\delta,\rho\in (0,1)$, $d\ge \delta n$ and $\rho<\rho_N(\delta)$, then}
$$ \lim_{n\to\infty} \bP\left(f_k(\Pi T^{n-1})=\binom{n}{k+1}\mbox{ for }k=0,\dots,\lfloor\rho d\rfloor\right)=1.$$

\vspace{2mm}

Thus, the probability that $\Pi T^{n-1}$ is $(k+1)$-neighborly for $k=1,\dots\lfloor\rho d\rfloor$, tends to one as $n\to\infty$. For further results in this direction, and hints to possible applications, we refer to Donoho and Tanner \cite{DT09}.

Somewhat earlier, Vershik and Sporyshev \cite{VS92} had considered the same model of random polytopes and investigated a weaker notion of neighborliness, roughly asking whether $f_k(\Pi T^{n-1})\ge \binom{n}{k+1}(1-\varepsilon)$ with $\varepsilon>0$, under a linearly coordinated growth of $d,k,n$.

We mention that a result of Baryshnikov and Vitale \cite{BV94} allows to transform the previous results into results about convex hulls of independent Gaussian samples. 

In the following, we want to take up Gale's original suggestion and study neighborliness properties of high-dimensional  polytopes when their Gale diagrams are generated at random. We recall briefly (a few more details will be given in Section \ref{sec2}) that with a sequence $(a_1,\dots,a_N)$ of points in $\R^d$, which affinely span $\R^d$, one can associate a Gale transform, which is a sequence $(\overline a_1,\dots,\overline a_N)$ of vectors in $\R^{N-d-1}$, linearly spanning $\R^{N-d-1}$. These vectors satisfy $\overline a_1+\dots+\overline a_N=o$ (where $o$ denotes the origin of $\R^{N-d-1}$). If $a_1,\dots,a_N$ are in affinely general position, then $\overline a_1,\dots,\overline a_N$ are in linearly general position. Any sequence $(\lambda_1\overline a_1,\dots,\lambda_N\overline a_N)$ with $\lambda_i>0$ for $i=1,\dots,N$ is a Gale diagram of $(a_1,\dots,a_N)$ (for the general definition of a Gale diagram, which is not needed here, we refer to \cite{McMS71}). Conversely, a sequence $(\overline a_1,\dots,\overline a_N)$ of vectors in $\R^{N-d-1}$, positively spanning the space, is the Gale diagram of a sequence $(a_1,\dots,a_N)$ in $\R^d$, in fact of many such sequences, but their convex hulls are combinatorially equivalent polytopes.

Now we assume that $\phi$ is a probability measure on $\R^{N-d-1}$ with the following properties:\\[1mm]
(a) $\phi$ is even (invariant under reflection in $o$),\\[1mm]
(b) $\phi$ assigns measure zero to each hyperplane through $o$.

\vspace{2mm}

\noindent{\bf Definition.} A $(\phi,N)$ {\em random Gale diagram} is a sequence $(X_1,\dots,X_N)$ of independent random vectors in $\R^{N-d-1}$ with distribution $\phi$, under the condition that
$$ o\in {\rm conv}\{X_1,\dots,X_N\}.$$

%\vspace{1mm}

Let $(X_1,\dots,X_N)$ be a realization of a $(\phi,N)$ random Gale diagram. Since $o\in {\rm conv}\{X_1,\dots,X_N\}$
and almost surely $X_1,\dots,X_N$ are in linearly general position (as follows from property (b) of the distribution $\phi$), we even have $o\in {\rm int\,conv}\{X_1,\dots,X_N\}$ a.s. Then (see, e.g., \cite[Thm. 1.1.14]{Sch14}), there are positive numbers $\lambda_i$ such that the sequence $(\overline a_1,\dots,\overline a_N)=(\lambda_1 X_1,\dots,\lambda_N X_N)$ satisfies $\overline a_1+\dots+\overline a_N=o$. Therefore, $(\overline a_1,\dots,\overline a_N)$ is the Gale transform of a sequence $(a_1,\dots,a_N)$ of points in $\R^d$. Let 
$$G_{d,N}:= {\rm conv}\{a_1,\dots,a_N\}.$$ 
The sequence $(a_1,\dots,a_N)$ and hence the polytope $G_{d,N}$ are not uniquely determined by $(\overline a_1,\dots,\overline a_N)$. However, it is determined by the Gale transform $(\overline a_1,\dots,\overline a_N)$, and in fact already by $(X_1,\dots,X_N)$, which points of $a_1,\dots,a_N$ are vertices of a face of $G_{d,N}$. Therefore, all polytopes $G_{d,N}$ which are determined by a given sequence $(X_1,\dots,X_N)$ are combinatorially isomorphic. In other words, we do not define random polytopes here, but random sets of combinatorially equivalent polytopes. We need not care about a measurable selection, since we are only interested in $f_k(G_{d,N})$, the number of $k$-faces of $G_{d,N}$, and this does not depend on the choice of $(a_1,\dots,a_N)$ with Gale transform $(\overline a_1,\dots,\overline a_N)$, but is determined by the sequence $(X_1,\dots,X_N)$. For that reason, $f_k(G_{d,N})$ is a well-defined random variable.

Let $k\in\{1,\dots,d-1\}$. Since $G_{d,N}$ is simplicial (a.s.), it is $(k+1)$-neighborly if and only if $f_k(G_{d,N})=\binom{N}{k+1}$, and if this holds, then $G_{d,N}$ is $j$-neighborly for $j\in\{2,\dots,k+1\}$. Therefore, we need only consider the random variable $f_k(G_{d,N})/\binom{N}{k+1}$. For this, we have the following results. 

Our first theorem is motivated by Theorem 1.3 of Donoho and Tanner \cite{DT10} (who consider a class of random cones), and it exhibits the same threshold, which we briefly recall. One defines
$$ H(x):= -x\log x-(1-x)\log(1-x), \quad 0\le x\le 1,$$
(with $0\log 0:=1$) and
$$ G(\delta,\rho):= H(\delta)+\delta H(\rho)-(1-\delta\rho)\log 2,\quad 0\le\delta,\rho\le 1.$$
If $\delta>1/2$, the function $G(\delta,\cdot)$ has a unique zero in $(0,1)$ (see also \cite[Lem. 6]{HS20}), which is denoted by $\rho_S(\delta)$. The graph of the function $\rho_S$ is depicted in \cite[Fig.1]{DT10}.

\begin{theorem}\label{T1.1}
Let $1/2<\delta<1$ and $0<\rho<1$ be given. Let $k<d<N-1$ be integers satisfying
$$ \frac{d}{N}\to\delta, \qquad \frac{k}{d}\to\rho\qquad\mbox{as }d\to\infty.$$
Then
$$ \lim_{d\to\infty}\bP\left(f_k(G_{d,N})=\binom{N}{k+1}\right)=1\quad\mbox{if } \rho<\rho_S(\delta).$$
\end{theorem}

If we ask only for the expectation (denoted by $\bE$) of $f_k(G_{d,N})$, we obtain a phase transition as in \cite{DT10}, with the same weak threshold, defined by 
$$ \rho_W(\delta):= \max\{0,2-\delta^{-1}\}, \qquad0<\delta<1.$$

\begin{theorem}\label{T1.2}
Let $0<\delta,\rho<1$ be given. Let $k<d<N-1$ be integers satisfying
$$ \frac{d}{N}\to\delta, \qquad \frac{k}{d}\to\rho\qquad\mbox{as }d\to\infty.$$
Then
$$ \lim_{d\to\infty}\frac{\bE f_k(G_{d,N})}{\binom{N}{k+1}}=\left\{\begin{array}{ll} 1 & \mbox{if }\rho<\rho_W(\delta),\\ 0 & \mbox{if }\rho>\rho_W(\delta).\end{array}\right.$$
\end{theorem}

After recalling some facts about Gale transforms in the next section, we show in Section \ref{sec3} how these theorems follow from asymptotic results obtained previously in a different context.

\section{Gale transforms}\label{sec2}

We first recall the essential facts about Gale transforms; more information can be found in the books by Gr\"unbaum \cite{Gru03}, McMullen and Shephard \cite{McMS71}, or Matou\v{s}ek \cite{Mat02}. We refer also to the survey given by McMullen \cite{McM79}. Let $(a_1,\dots,a_N)$ be a sequence of points in $\R^d$, which affinely span $\R^d$. We write vectors as ordered tuples of coordinates with respect to an orthonormal basis, for example,
$$ \begin{array}{lll} a_1 &=& (\alpha_{11}, \alpha_{21},\dots,\alpha_{d1}),\\
\vdots && \\
 a_N &=& (\alpha_{1N}, \alpha_{2N},\dots,\alpha_{dN}).\end{array} $$
Using theses vectors, we form the $N\times N$ matrix
$$ M= \left(\begin{array}{llll}
\alpha_{11} & \alpha_{12} & \dots & \alpha_{1N}\\
\vdots & \vdots & &\vdots \\
\alpha_{d1} & \alpha_{d2} & \dots & \alpha_{dN}\\
1 & 1 & \dots & 1\\
\beta_{11} & \beta_{12} & \dots & \beta_{1N}\\
\vdots & \vdots & &\vdots \\
\beta_{N-d-1,1} & \beta_{N-d-1,2} & \dots & \beta_{N-d-1,N}\end{array}\right)$$
in such a way that it has rank $N$ and the last $N-d-1$ rows are orthogonal to the first $d+1$ rows. Then $(\overline a_1,\dots,\overline a_N)$ with
$$ \begin{array}{lll} \overline a_1 &=& (\beta_{11}, \beta_{21},\dots,\beta_{N-d-1,1}),\\
\vdots && \\
 \overline a_N &=& (\beta_{1N}, \beta_{2N},\dots,\beta_{N-d-1,N}),\end{array} $$
is a Gale transform of $(a_1,\dots,a_N)$. In other words, the vectors $a_1,\dots,a_N$ are the columns of the upper $d\times N$ submatrix of $M$. The affine dependences $(\lambda_1,\dots,\lambda_N)$ of $a_1,\dots,a_N$, defined by 
$$ \sum_{i=1}^N \lambda_ia_i=o,\qquad  \sum_{i=1}^N \lambda_i=0,$$
form a vector space of dimension $N-d-1$. A basis of this vector space makes up the last $N-d-1$ rows of $M$, and $\overline a_1,\dots,\overline a_N$ are the columns of the lower $(N-d-1)\times N$ submatrix of $M$.

It is clear from this description that any sequence $(\overline a_1,\dots,\overline a_N)$ of vectors in $\R^{N-d-1}$, which linearly span $\R^{N-d-1}$ and satisfy
$$ \overline a_1+\dots+\overline a_N=o,$$
is the Gale transform of a sequence $(a_1,\dots,a_N)$ of points in $\R^d$. Of course, the latter sequence is not uniquely determined, as there is freedom in the choice of bases.

The points $a_1,\dots,a_N$ are in affinely general  position in $\R^d$ if and only if the vectors $\overline a_1,\dots,\overline a_N$ are in linearly general position in $\R^{N-d-1}$. For vectors $x_1,\dots,x_m$ in linearly general position, the relations $o\in {\rm relint\,conv}\{x_1,\dots,x_m\}$ and $o\in {\rm conv}\{x_1,\dots,x_m\}$ are equivalent. Therefore, we can state the following crucial lemma, for which we refer, e.g., to \cite{Gru03} or \cite{Mat02}. 

\begin{lemma}\label{L2.1}
Let $(\overline a_1,\dots,\overline a_N)$ be a Gale transform of the sequence $(a_1,\dots,a_N)$ of points in $\R^d$ ($N\ge d+1$), where $\overline a_1,\dots,\overline a_N$ linearly span $\R^{N-d-1}$ and are in linearly general position. Let $P={\rm conv}\{a_1,\dots,a_N\}$. Then $P$ is a simplicial polytope, and for $k\in\{0,\dots,d-1\}$ and any $(k+1)$-element subset $I\subset\{1,\dots,N\}$ we have:
$$ {\rm conv}\{a_i:i\in I\} \mbox{ is a $k$-face of $P$} \,\Leftrightarrow\, o\in {\rm conv}\{\overline a_j:j\notin I\}.$$
\end{lemma}

\section{Proofs of the theorems}\label{sec3}

Let integers $d\ge 1$ and $N\ge d+1$ be given. In the following, we assume that $\phi$ is a distribution on $\R^{N-d-1}$ as specified above, and that $(X_1,\dots,X_N)$ is a $(\phi,N)$ {\em random Gale diagram}. A given realization satisfies $o\in {\rm conv}\{X_1,\dots,X_N\}$, and since it is almost surely in linearly general position, we have even $o\in {\rm int\,conv}\{X_1,\dots,X_N\}$. Consequently, there are positive numbers $\lambda_1,\dots,\lambda_N$ such that the sequence $(\overline a_1,\dots,\overline a_N):= (\lambda_1 X_1,\dots,\lambda_N X_N)$ satifies $\overline a_1+\dots+\overline a_N=o$. Therefore, it is a Gale transform of a sequence $(a_1,\dots,a_N)$ of points in $\R^d$. Let $G_{d,N}={\rm conv}\{a_1,\dots,a_N)$. As mentioned, $G_{d,N}$ is a simplicial polytope. For $k\in\{0,\dots,d-1\}$, let $f_k(G_{d,N})$ denote the number of its $k$-faces. Since 
$$ o\in {\rm conv}\{\overline a_j:j\notin I\} \,\Leftrightarrow\, o\in {\rm conv}\{X_j:j\notin I\},$$
it follows from Lemma \ref{L2.1} that
$$ f_k(G_{d,N}) = \sum_{1\le i_1<\dots<i_{N-k-1}\le N} {\mathbbm 1}\left\{o\in{\rm conv}\{X_{i_1},\dots,X_{i_{N-k-1}}\}\right\}.$$
We conclude that $f_k(G_{d,N})$ is a random variable, depending only on $(X_1,\dots,X_N)$, and that its expectation is given by a sum of conditional probabilities,
\begin{eqnarray}\label{3.1} 
&& \bE f_k(G_{d,N})\\
&&= \sum_{1\le i_1<\dots<i_{N-k-1}\le N} \bP\left(o\in{\rm conv}\{Y_{i_1},\dots,Y_{i_{N-k-1}}\}\mid o\in{\rm conv}\{Y_1,\dots,Y_N\}\right),\nonumber
\end{eqnarray}
where $Y_1,\dots,Y_N$ are independent random vectors in $\R^{N-d-1}$ with distribution $\phi$. Since
$$ o\in{\rm conv}\{Y_{i_1},\dots,Y_{i_{N-k-1}}\}\,\Rightarrow\, o\in{\rm conv}\{Y_1,\dots,Y_N\},$$
we obtain
$$ \bE f_k(G_{d,N})= \frac{\sum_{1\le i_1<\dots<i_{N-k-1}\le N} \bP\left(o\in{\rm conv}\{Y_{i_1},\dots,Y_{i_{N-k-1}}\}\right)}
{\bP\left(o\in{\rm conv}\{Y_1,\dots,Y_N\}\right)}.$$

The probabilities occurring here can be determined with the aid of the following lemma. It is due to Wendel \cite{Wen62}. The proof is reproduced in \cite[Thm. 8.2.1]{SW08}.

\begin{lemma}\label{L3.1}
If $Y_1,\dots,Y_M$ are i.i.d. random vectors in $\R^r$ with a symmetric distribution which is zero on hyperplanes through $o$, then
$$ P_{r,M}:= \bP\left(o\notin{\rm conv}\{Y_1,\dots,Y_M\}\right) = \frac{1}{2^{M-1}} \sum_{i=0}^{r-1}\binom{M-1}{i}.$$
\end{lemma}

Since $P_{r,M}+P_{M-r,M}=1$, we have
$$ \bP\left(o\in{\rm conv}\{Y_1,\dots,Y_M\}\right) = P_{M-r,M}.$$
This gives
\begin{equation}\label{3.2} 
\bE f_k(G_{d,N})=  \binom{N}{k+1} \frac{P_{d-k,N-k-1}}{P_{d+1,N}}.
\end{equation}

We can now take advantage of the fact that the expressions appearing in (\ref{3.2}) have shown up in a different situation, and their asymptotic behavior has already been investigated. Let $\psi$ be a probability distribution on $\R^{d+1}$, with properties corresponding to those that $\phi$ has on $\R^{N-d-1}$. For $N\in\N$, the $(\psi,N)$ Cover--Efron cone $C_{d+1,N}$ is a random cone, defined as the positive hull of $N$ independent random vectors $Z_1,\dots,Z_N$ in $\R^{d+1}$ with distribution $\psi$, under the condition that ${\rm pos}\{Z_1,\dots,Z_N\}\not=\R^{d+1}$. If $f_{k+1}(C)$ denotes the number of $(k+1)$-dimensional faces of a polyhedral cone $C$, we have 
\begin{equation}\label{3.3}
\frac{\bE f_{k+1}(C_{d+1,N})}{\binom{N}{k+1}} = \frac{P_{d-k,N-k-1}}{P_{d+1,N}}.
\end{equation}
This was proved in \cite[(3.3)]{CE67}, and also in \cite[(27)]{HS20}. From (\ref{3.2}) and (\ref{3.3}) it follows that
\begin{equation}\label{3.4} 
\bE f_k(G_{d,N}) = \bE f_{k+1}(C_{d+1,N}).
\end{equation}
Thus, results on expected face numbers of Cover--Efron cones immediately imply results on expected face numbers of the polytopes $G_{d,N}$.

\vspace{2mm}

\noindent{\bf Remark.} Instead of viewing the sequence $(X_1,\dots,X_N)$ as a Gale diagram of a sequence $(a_1,\dots,a_N)$ in $\R^d$, we can also view it as a linear transform of a sequence $(b_1,\dots,b_N)$ in $\R^{d+1}$. (For linear transforms, we refer to Shephard \cite{She71}, where they are called linear representations, and to McMullen \cite{McM79}.) Define
$$ D_{d+1,N} := {\rm pos}\{b_1,\dots,b_N\}.$$
Then $D_{d+1,N}$ is a polyhedral cone, different from $\R^{d+1}$ (by \cite{She71}, Corollary to Theorem 1). Again, $D_{d+1,N}$ is not uniquely determined by the realization $X_1,\dots,X_N$, but its number of $(k+1)$-faces, $f_{k+1}(D_{d+1,N})$, is determined by  $(X_1,\dots,X_N)$ and thus is a well-defined random variable. One finds that
$$ \frac{\bE f_{k+1}(D_{d+1,N})}{\binom{N}{k+1}} = \frac{P_{d-k,N-k-1}}{P_{d+1,N}} = \frac{\bE f_{k+1}(C_{d+1,N})}{\binom{N}{k+1}}.$$
Thus, from the viewpoint of expected face numbers and neighborliness properties, this model of random cones is not different from the model of Cover--Efron cones.

\vspace{2mm}

\noindent{\em Proof of Theorem} \ref{T1.1}

As above, we assume that $Y_1,\dots,Y_N$ are independent random vectors in $\R^{N-d-1}$ with distribution $\phi$. Let ${\mathcal I}_k$ be the set of $(N-k-1)$-element subsets of $\{1,\dots,N\}$. The conditional probability
$$ p(I):= \bP\left(o\notin {\rm conv}\{Y_i:i\in I\}\mid o\in {\rm conv}\{Y_1,\dots,Y_N\}\right)$$
is independent of the choice of $I\in{\mathcal I}_k$ and can therefore be denoted by $p(I_0)$, for a fixed $I_0\in{\mathcal I}_k$. By Boole's inequality,
\begin{eqnarray*}
&& \bP\left(o\notin {\rm conv}\{Y_i:i\in I\}\mbox{ for some $I\in{\mathcal I}_k$}\mid o\in {\rm conv}\{Y_1,\dots,Y_N\}\right)\\
&&\le \sum_{I\in{\mathcal I}_k} \bP\left(o\notin {\rm conv}\{Y_i:i\in I\} \mid o\in {\rm conv}\{Y_1,\dots,Y_N\}\right)\\
&&=\binom{N}{k+1}p(I_0).
\end{eqnarray*}
It follows that
\begin{eqnarray*}
&& \bP\left(f_k(G_{d,N})=\binom{N}{k+1}\right)\\
&&= 1-\bP\left(o\notin {\rm conv}\{Y_i:i\in I\}\mbox{ for some $I\in{\mathcal I}_k$}\mid o\in {\rm conv}\{Y_1,\dots,Y_N\}\right)\\
&& \ge 1-\binom{N}{k+1}p(I_0).
\end{eqnarray*}
Since
$$ \bE f_k(G_{d,N})=\binom{N}{k+1}\bP\left(o\in {\rm conv}\{Y_i:i\in I_0\} \mid o\in{\rm conv}\{Y_1,\dots,Y_N\}\right)$$
by (\ref{3.1}), we have
$$ p(I_0) =1-\frac{\bE f_k(G_{d,N})}{\binom{N}{k+1}}.$$
Now (\ref{3.4}) allows us to write
$$ p(I_0)=1-\frac{\bE f_{k+1}(C_{d+1,N})}{\binom{N}{k+1}} =\frac{A}{1+A}\le A,$$
with $A$ defined by \cite[(13)]{HS20}, and here evaluated at $d+1,N,k+1$. The proof is now completed precisely as that of Theorem 6 in \cite{HS20}, yielding that $p(I_0)\to 0$ as $d\to\infty$. Note that the change from $(d,k)$ to $(d+1,k+1)$ does not alter the assumptions of \cite[Thm. 6]{HS20}. \hfill $\Box$

\noindent{\em Proof of Theorem} \ref{T1.2}.
In view of (\ref{3.4}), the proof is the same as that for Theorem 5 in \cite{HS20}. \hfill $\Box$

\noindent Author's address:\\[2mm]
Rolf Schneider\\Mathematisches Institut, Albert-Ludwigs-Universit{\"a}t\\D-79104 Freiburg i.~Br., Germany\\E-mail: rolf.schneider@math.uni-freiburg.de

\end{document}